\documentclass[11pt,reqno]{amsart}

\usepackage{amsmath,amssymb,amsthm,mathrsfs,graphicx,xcolor}

\theoremstyle{plain}
\newtheorem{theorem}{Theorem}
\newtheorem{lemma}[theorem]{Lemma}

\newtheorem{question}[theorem]{Question}

\theoremstyle{definition}
\newtheorem*{remark}{Remark}

\begin{document}

\title[On a generalization of Busemann's intersection inequality]
      {On a generalization of Busemann's intersection inequality}

\author[V.~Yaskin]{Vlad Yaskin}
\address{Vladyslav Yaskin, Department of Mathematical and Statistical Sciences, University of Alberta, Edmonton, Alberta, T6G 2G1, Canada}
\email{yaskin@ualberta.ca}

\subjclass[2010]{52A30 (primary), and 42B10 (secondary)}

\thanks{The author was partially supported by NSERC}

\keywords{Intersection body, Busemann's intersection inequality, Fourier transform, spherical harmonics, interpolation of operators}

\begin{abstract} Busemann's intersection inequality gives an upper bound for the volume of the intersection body of a star body in terms of the volume of the body itself. Koldobsky, Paouris, and Zymonopoulou asked if there is a similar result for $k$-intersection bodies. 
We solve this problem for star bodies that are close to the Euclidean ball in  the Banach-Mazur distance. We also improve a bound obtained by Koldobsky, Paouris, and Zymonopoulou for general star bodies in the case when $k$ is proportional to the dimension.
\end{abstract}

\maketitle

\section{Introduction}
 We say that a set $K$ in $\mathbb{R}^n$ is {\it star-shaped}  if for every $x\in K$ the closed line segment connecting $x$ to the origin lies in $K$. A compact set $K$ in $\mathbb{R}^n$ is called a {\it star body}  if it is star-shaped and its {\it radial function}  defined by
\begin{align*}
\rho_K(\xi) = \max\{a \ge 0 :a \xi \in K\}, \quad \xi\in S^{n-1},
\end{align*}
is positive and continuous. Geometrically, $\rho_K(\xi)$ is the distance from the origin to the point on the boundary in the direction of $\xi$.

Let $K$ and $L$ be  star bodies in $\mathbb R^n$. Following Lutwak \cite{L}, we say that $L$ is the {\it intersection body} of $K$ if $$\rho_{L}(\xi)  = \mathrm{vol}_{n-1} ( K\cap\xi^\perp ),$$
for every $\xi\in S^{n-1}$.

For any star body $K$ its intersection body always exists and  we will denote it  by $I(K)$. The well-known
Busemann intersection inequality asserts that $$\mathrm{vol}_{n}(I(K))\le \frac{\kappa_{n-1}^n}{\kappa_n^{n-2}} \mathrm{vol}_{n}(K)^{n-1}$$
with equality if and only if $K$ is a centered ellipsoid; see, e.g., \cite[p.~373]{Ga}. Here and below, $$\kappa_p = \frac{\pi^\frac{p}{2} }{\Gamma\left( 1 + \frac{p}{2} \right)},$$ which equals the volume of the unit Euclidean ball in $\mathbb R^p$ when $p$ is a positive integer.

Koldobsky introduced a generalization of the notion of intersection body; see \cite[p.~75]{K}. Let $K$ and $L$ be origin-symmetric star bodies in $\mathbb R^n$ and $k$ be an integer, $1\le k\le n-1$. We say that $L$ is the $k$-intersection body of $K$ if $$\mathrm{vol}_{k}(L\cap H) = \mathrm{vol}_{n-k}(K \cap H^\perp)$$ for every $k$-dimensional subspace $H$ of $\mathbb R^n$.

For a given origin-symmetric star body $K$ its $k$-intersection body may not exist, but when it does, we will denote it by $I_k (K)$.

Koldobsky, Paouris, and Zymonopoulou \cite{KPZ} asked whether an analogue of Busemann's intersection inequality holds for $k$-intersection bodies. Namely, if $K$ is an origin-symmetric star body in $\mathbb R^n$ whose $k$-intersection body exists and such that $\mathrm{vol}_{n}(K)=\mathrm{vol}_{n}( B_2^n)$, is it true that
\begin{equation}\label{KPZ-conj}
\mathrm{vol}_{n}(I_k (K))\le \mathrm{vol}_{n}(I_k(B_2^n))?
\end{equation}
They proved that 
\begin{equation}\label{KPZ} \left(\frac{\mathrm{vol}_{n}(I_k (K))}{ \mathrm{vol}_{n}(I_k(B_2^n))}\right)^{1/n} \le c \min\{\log n, k\log k\},
\end{equation}
for some absolute constant $c$, under the assumption that $I_k (K)$ is a convex body. 

We will slightly modify and extend their conjecture.
First of all, let us write (\ref{KPZ-conj}) somewhat differently. If $L$ is the $k$-intersection body of $K$, and we do not put any restrictions on the volume of $K$, then (\ref{KPZ-conj}) is equivalent to
\begin{equation}\label{KPZ-conj1}
\left(\mathrm{vol}_{n}(L)\right)^k \le C_{n,k} \left(\mathrm{vol}_{n}(K)\right)^{n-k},
\end{equation}
where $C_{n,k}$ is an appropriate constant so that the latter inequality becomes equality in the case of centered balls. 

It follows directly from the definition that if $L$ is the $k$-intersection body of $K$, then $K$ is the $(n-k)$-intersection body of $L$. This means that if inequality (\ref{KPZ-conj1}) is true for $0<k<n/2$, then the reversed inequality should be true for $n/2<k<n$.

Further, in terms of the Fourier transform the condition that $L$ is the $k$-intersection body $K$ can be written as follows: 
\begin{equation}\label{FT-k-int}
\|\theta\|_L^{-k} = \frac{k}{(2\pi)^k(n-k)} \left( \|\cdot\|_K^{-n+k} \right)^\wedge(\theta), \quad \theta \in S^{n-1};
\end{equation}
  see \cite[Theorem 4.6]{K}. Here $\|\cdot\|_K$ denotes the {\it Minkowski functional} of $K$ and is defined by
\begin{align*}
\|x\|_K=\min \{a\ge 0: x \in aK \}, \qquad x\in \mathbb R^n.
\end{align*}

Note that   relation (\ref{FT-k-int}) allows us to write (\ref{KPZ-conj1}) as follows:
 \begin{equation}\label{FT-KPZ}  \int_{S^{n-1}} \left[ (\|\cdot\|_K^{-n+k})^\wedge(\theta) \right]^{n/k} d\theta \le c_{n,k} \left(\mathrm{vol}_{n}(K)\right)^{\frac{n}{k}-1},
\end{equation}
where $c_{n,k}$ is a constant that turns (\ref{FT-KPZ}) into equality for centered balls. The exact value of this constant will be computed later.

 As one can see,   (\ref{FT-KPZ}) makes sense even if $k$ is not an integer. Thus we will write (\ref{FT-KPZ}) for a larger set of values by using  a real number  $p$ instead of $k$. Let us also note that the assumption that the Fourier transform of $\|\cdot\|_K^{-n+k}$ is a positive continuous function is very restrictive. So this assumption will be dropped.

Summarizing all of the above remarks, let us now state the  question in the following form.

\begin{question}\label{Conj1} Let $0<p<n$ be a real number and $K$ be an origin-symmetric star body in  $\mathbb R^n$.
Are the following inequalities true? 
\begin{equation}\label{C1}  \int_{S^{n-1}} \left| (\|\cdot\|_K^{-n+p})^\wedge(\theta) \right|^{n/p} d\theta \le c_{n,p} \left(\mathrm{vol}_{n}(K)\right)^{\frac{n}{p}-1}, \quad \mbox{if } \, p< n/2,
\end{equation}
and
\begin{equation}\label{C2}  \int_{S^{n-1}} \left| (\|\cdot\|_K^{-n+p})^\wedge(\theta) \right|^{n/p} d\theta \ge c_{n,p} \left(\mathrm{vol}_{n}(K)\right)^{\frac{n}{p}-1}, \quad \mbox{if } \, p> n/2,
\end{equation}
with equality if and only if $K$ is a centered ellipsoid. 
\end{question}

The case $p=n/2$ is omitted above since these inequalities become equalities for every origin-symmetric star body $K$. This is just  an application of the spherical version of Parseval's formula (see \cite[p.66]{K}):
$$   \int_{S^{n-1}} \left[ (\|\cdot\|_K^{-n/2})^\wedge(\theta) \right]^{2} d\theta = (2\pi)^n \int_{S^{n-1}}  \|\theta\|_K^{-n }   d\theta = n (2\pi)^n  \mathrm{vol}_{n}(K). $$

It is interesting to note that the  conjectured inequalities (\ref{C1}) and (\ref{C2}) have connections to other known inequalities. 
First of all, let us repeat that the case $p=1$ corresponds to the Busemann intersection inequality. Thus Question \ref{Conj1} has an affirmative answer for $p=1$ and $p=n-1$.
Let us now look at the case when $p<1$ and $p$ is not an even integer. For these values of $p$ the Fourier transform of $\|\cdot\|_K^{-n+p}$ can be expressed as follows:
\begin{equation}\label{IntFT}(\|\cdot\|_K^{-n+p})^\wedge(\theta) = \frac{\pi (n-p)}{2\Gamma(1-p) \sin(\pi p/2) } \int_K |\langle x,\theta \rangle|^{-p}\, dx;
\end{equation}
see \cite[Corollary 3.15]{K}.

The reader may recognize such integrals: they appear, for example, in the definition of polar $q$-centroid bodies (see \cite{LZ} for $q\ge 1$ and \cite{YY} for $-1<q<1$). Let $K\subset\mathbb{R}^n$ be a star body and $q>-1$, $q\neq 0$. The polar $q$-centroid body of $K$ is the star body $\Gamma_q^* K$ given by
\begin{align*}
\|x\|_{\Gamma_q^* K} = \left( \frac{1}{ \mbox{vol}_n(K)} \int_K |\langle x,y\rangle|^q\, dy \right)^\frac{1}{q} , \quad x\in\mathbb{R}^n\setminus\{0\}.
\end{align*}
Note that in \cite{LZ} the normalization is different.

Lutwak and Zhang \cite{LZ} have shown that if $K\subset\mathbb{R}^n$ is a star body and $q\ge 1$, then
\begin{align}\label{LZ}
\mbox{vol}_n(K) \mbox{vol}_n(\Gamma_q^*K) \leq \mbox{vol}_n(B_2^n) \mbox{vol}_n(\Gamma_q^* B_2^n),
\end{align}
with equality if and only if $K$ is an   ellipsoid centered at the origin.

By virtue of formula (\ref{IntFT}) one can check that, in the case when $q$ is not an even integer, the  Lutwak-Zhang inequality is equivalent to
$$  \int_{S^{n-1}} \left| (\|\cdot\|_K^{-n-q})^\wedge(\theta) \right|^{-n/q} d\theta \le |c_{n,-q}| \left(\mathrm{vol}_{n}(K)\right)^{-\frac{n}{q}-1} . 
$$
 Thus, Question \ref{Conj1} can also be viewed as an extension of the Lutwak-Zhang inequality. 
 
In this paper  we   show that Question \ref{Conj1} has an affirmative answer when the body $K$ is sufficiently close to the Euclidean ball in the Banach-Mazur distance. 
For general star bodies we will obtain an improvement of   inequality (\ref{KPZ}) when $k>cn/\log^2n$.

\section{Preliminaries}\label{Prelim}

Two of the main tools used in this paper are the Fourier transform  of distributions and spherical harmonics. The reader is referred to the books \cite{K} and \cite{Gr}  for  detailed discussions of such techniques. We will just briefly mention some important facts. Let $f$ be a continuous function on the sphere $S^{n-1}$ and   consider its homogeneous extension to $\mathbb R^n\setminus\{0\}$ of degree $-n+p$, where $0<p<n$.  
We can think of $|x|_2^{-n+p} f(x/|x|_2)$ as a distribution acting on test functions by integration and therefore we can define its Fourier transform  in the distributional sense.  If $f\in C^\infty(S^{n-1})$, then the Fourier transform of $|x|_2^{-n+p} f(x/|x|_2)$ is equal to a homogeneous of degree $-p$ function, that is infinitely smooth on $\mathbb R^n\setminus\{0\}$; see \cite[Lemma  3.16]{K}.  
Thus we can introduce a linear operator $I_p: C^\infty(S^{n-1})\to  C^\infty(S^{n-1})$, that maps  a function $f$ to the function $I_pf$
equal to the restriction to the sphere of the Fourier transform  of $\frac{\Gamma\left( \frac{ n-p}{2}\right)}{2^p\pi^\frac{n}{2} \Gamma\left(\frac{ p}{2}\right)} |x|_2^{-n+p} f(x/|x|_2)$. The coefficient in front of the latter function is chosen in such a way that $I_p(1)=1$, as will be shown later.

When $0<p<1$,  $I_p$ has the following integral representation, which is just  (\ref{IntFT}) with an appropriate normalization.
\begin{align}\label{IntFT1}I_pf(\theta)& = \frac{\Gamma\left( \frac{ n-p}{2}\right)}{2^p\pi^\frac{n}{2} \Gamma\left(\frac{ p}{2}\right)} \frac{\pi  }{2\Gamma(1-p) \sin(\pi p/2) } \int_{S^{n-1}} |\langle x,\theta \rangle|^{-p} f(x) \, dx\nonumber\\
& = \frac{ \sqrt{\pi} \Gamma\left( \frac{ n-p}{2}\right)}{ \Gamma\left(\frac{ n }{2}\right) \Gamma\left(\frac{ 1-p}{2}\right)}   \int_{S^{n-1}} |\langle x,\theta \rangle|^{-p} f(x) \, d\sigma(x),
\end{align}
where $\sigma$ is the rotationally invariant probability measure on the sphere. To compute the coefficient in front of the above integral we used \cite[Lemma 2.18]{K}.

For a function $f$ on the sphere, let $\sum_{m=0}^\infty H_m$ denote   its  spherical harmonic expansion, where each $H_m$ is a spherical harmonic of degree $m$ in dimension $n$.  
 If $f$ is   even, then it has only harmonics of even degrees in its expansion:  $\sum_{m\ge 0, \, m \, \mathrm{even}}  H_m$. Furthermore,  $I_pf$ is also even and its spherical harmonic expansion is given by
 $$\sum_{m\ge 0,\, \atop  m \, \mathrm{even}} \lambda_m(n,p) H_m,$$
 where 
$$\lambda_m(n,p) = \frac{\Gamma\left( \frac{ n-p}{2}\right)\Gamma\left(\frac{m+p}{2}\right)}{\Gamma\left(\frac{ p}{2}\right)\Gamma\left( \frac{m+n-p}{2}\right)};$$
see, e.g., \cite{GYY} (note that the normalization of $I_p$ is different there).

Note that $\lambda_0(n,p)=1$ for all $p\in (0,n)$, and thus  $I_p(1)=1$. Another important observation is that $|\lambda_m(n,p)|>|\lambda_{m+2} (n,p)|$ for all even $m\ge 0$ when $0<p<n/2$ (and the inequality gets reversed  when $n/2<p<n$). In particular, for $0<p<n/2$, we have  $$\|I_pf\|_2^2 =  \sum_{m\ge 0,\,  \atop  m \, \mathrm{even}}  \lambda_m^2(n,p) \|H_m\|_2 \le   \| f\|_2^2 .$$
Hence, $I_p$  (when $0<p<n/2$) is well defined as a 
linear operator from $L_{even}^2(S^{n-1})$ to $L_{even}^2(S^{n-1})$ with the operator norm  equal to 1.

Here and below we denote by $\|f\|_q$ the $L^q(S^{n-1})$-norm of $f$:
$$\|f\|_q = \left(\int_{S^{n-1}} |f(\theta)|^q d\sigma(\theta)\right)^{1/q} .$$

Let us finally remark that in terms of the operator $I_p$ the conjectured inequality (\ref{C1})   can be written as follows.

\noindent {\bf Question 1 (Reformulated).} {\it	
Let $K$   be an origin-symmetric star body in $\mathbb R^n$ and   $0<p<n/2$. Is it true that  \begin{equation}\label{MainIneq} \int_{S^{n-1}} \left| I_p (\|\cdot\|_K^{-n+p}) (\theta) \right|^{n/p} d\sigma(\theta) \le     \left(\kappa_n\right)^{1-\frac{n}{p}} \left(\mathrm{vol}_n(K) \right)^{\frac{n}{p}-1},  
\end{equation}
with equality if and only if $K$ is a centered ellipsoid?}

Below we will not discuss inequality (\ref{C2}), since, at least in the language of $k$-intersection bodies, the cases $0<k<n/2$ and $n/2<k<n$ are   the same.

\section{Main Results}\label{MainRes}

Following the discussion in the previous section, we will now show that $I_p$ can be extended to a bounded linear operator on a larger space of functions.

\begin{theorem}\label{interp} Let $0<p<n/2$. 
	$I_p$ is a bounded linear operator from $L_{even}^{n/(n-p)}(S^{n-1})$ to $L_{even}^{n/p}(S^{n-1})$.
	Moreover, 
	$$\|I_{p}\|_{L_{even}^{n/(n-p)}(S^{n-1})\to L_{even}^{n/p}(S^{n-1})} \le  \left(\frac{n}{2p}\right)^{\frac{p}2}  . $$

	\end{theorem}
\proof

We will follow the Stein  interpolation theorem; see \cite{S} or \cite[Chapter V]{SW}. First observe that for every $f,g\in C^\infty(S^{n-1})$ the mapping  $$p\mapsto \int_{S^{n-1}} I_pf(x) g(x) d\sigma(x)$$ is analytic in the strip  $0< \Re(p)< n/2$, and   continuous up to the boundary of the strip. In addition, 
$$\left| \int_{S^{n-1}} I_pf(x) g(x) d\sigma(x)\right| \le \|I_pf\|_2  \|g\|_2 \le \|f\|_2  \|g\|_2.$$
We will now investigate the behavior of $I_p$ on the boundary of the strip.  
In order to be consistent with the notation in Stein's theorem we will write the strip in the form $0\le \Re(z)\le 1$, where $z=\frac2n p$. Thus, the boundary of the strip consists of two lines $z = is$ and  $z = 1+is$, where $s\in \mathbb R$.
Therefore, when $ \Re(z)=1$, i.e.,  $p= \frac{n}{2}(1+is)$,  we get
$$\lambda_m\left(n,{  \frac{n}{2}}(1+is)\right) = \frac{\Gamma\left( \frac{ n }{4}-\frac{ is }{4}\right)\Gamma\left( \frac{m }{2}+\frac{ n }{4}+\frac{ is }{4}\right)}{\Gamma\left( \frac{ n }{4}+\frac{ is }{4}\right)\Gamma\left( \frac{m }{2}+\frac{ n }{4}-\frac{ is }{4}\right)}.$$
Since $\Gamma(a+ib)$ is the complex conjugate of $\Gamma(a-ib)$, we see that \linebreak $|\lambda_m(n,  \frac{n}{2}(1+is))|=1$   for all real $s$, and so
$$\|I_{\frac{n}{2}(1+is)}f\|_2 = \| f\|_2.$$

Now consider the case $  \Re(z)=0$. Extending (\ref{IntFT1}) to $p= \frac{n}{2} is $, $s\in \mathbb R$, we get
 $$\|I_{\frac{n}{2}is}f\|_{\infty} \le  \frac{ \sqrt{\pi} |\Gamma\left( \frac{ n-isn/2}{2}\right)| }{\Gamma\left( \frac{ n}{2}\right) |\Gamma\left(\frac{ 1-isn/2}{2}\right)|}   \int_{S^{n-1}} | f(x) |\, d\sigma(x).$$

Let $A(s)$ be the  coefficient in front of the latter integral. To estimate $A(s)$ we will distinguish two cases according to the parity of $n$. If $n$ is odd, then 
\begin{align*} A(s)& =  \frac{ \sqrt{\pi} |\Gamma\left( \frac{ n-isn/2}{2}\right)| }{\Gamma\left( \frac{ n}{2}\right) |\Gamma\left(\frac{ 1-isn/2}{2}\right)|} \\
& =  \frac{\sqrt{\pi}\left|\frac{ n-2-isn/2}{2}\right|\cdot\left| \frac{ n-4-isn/2}{2}\right|\cdots \left|\frac{ 1-isn/2}{2}\right|\cdot  \left|\Gamma\left( \frac{ 1-isn/2}{2}\right)\right| } {   \frac{n-2}2 \cdot\frac{n-4}2\cdots\frac{1}2\sqrt{\pi}\cdot |\Gamma\left(\frac{ 1-isn/2}{2}\right)|}\\
&=\left|\frac{ n-2-isn/2}{n-2}\right|\cdot\left| \frac{ n-4-isn/2}{n-4}\right|\cdots \left|\frac{ 1-isn/2}{1}\right|\\
&= \prod_{k=0}^{(n-3)/2} \left(1+\frac{s^2 n^2}{4(2k+1)^2}\right)^{1/2}\le \prod_{k=0}^{\infty} \left(1+\frac{s^2 n^2}{4(2k+1)^2}\right)^{1/2} \\
&=  \left( \cosh \frac{\pi sn}{4} \right)^{1/2}\le {e^{ \frac{\pi |s| n}8}}.
\end{align*}
Above we used a representation of $\cosh$ as an infinite product. See e.g., \cite[VII, \S 5-6]{C} for details on the Weierstrass factorization theorem.

If $n$ is even, the argument is similar, but we will additionally need the following formulas:
$$|\Gamma(1+is)| = \left(\frac{\pi s}{\sinh \pi s}\right)^{1/2}, \qquad \left|\Gamma\left(\frac12+is\right)\right| = \left(\frac{\pi }{\cosh \pi s}\right)^{1/2},$$
which can be obtained from the Euler reflection formula; see \cite[p.~9 and p.~22]{AAR} for some details.
The above formulas imply 
\begin{equation*} \frac{|\Gamma\left( \frac{2-isn/2}{2}\right)|}{ |\Gamma\left(\frac{ 1-isn/2}{2}\right)|} = \left( \frac{sn}4 \coth \frac{\pi s n}4\right)^{1/2} .
\end{equation*}

Then we have, for even $n$, 
\begin{align*}A(s)& =  \frac{ \sqrt{\pi} |\Gamma\left( \frac{ n-isn/2}{2}\right)| }{\Gamma\left( \frac{ n}{2}\right) |\Gamma\left(\frac{ 1-isn/2}{2}\right)|} \\
& =  \frac{\sqrt{\pi} \left|\frac{ n-2-isn/2}{2}\right|\cdot\left| \frac{ n-4-isn/2}{2}\right|\cdots \left|\frac{ 2-isn/2}{2}\right|\cdot  \left|\Gamma\left( \frac{ 2-isn/2}{2}\right)\right|  } {    \frac{n-2}2 \cdot\frac{n-4}2\cdots\frac{2}2 \cdot |\Gamma\left(\frac{ 1-isn/2}{2}\right)| }\\ 
&=  \sqrt{\pi}\left|\frac{ n-2-isn/2}{n-2}\right|\cdot\left| \frac{ n-4-isn/2}{n-4}\right|\cdots \left|\frac{ 2-isn/2}{2}\right| \left( \frac{sn}4 \coth \frac{\pi s n}4\right)^{1/2}\\
&= \sqrt{\pi} \left( \frac{sn}4 \coth \frac{\pi s n}4\right)^{1/2} \prod_{k=1}^{(n-2)/2} \left(1+\frac{n^2 s^2}{4 (2k)^2}\right)^{1/2}\\
&\le \sqrt{\pi} \left( \frac{sn}4 \coth \frac{\pi s n}4\right)^{1/2} \prod_{k=1}^{\infty} \left(1+\frac{n^2 s^2}{4 (2k)^2}\right)^{1/2}\\
& = \sqrt{\pi}  \left( \frac{sn}4 \coth \frac{\pi s n}4\right)^{1/2}   \left(\frac{4\sinh \frac{\pi s n}4}{\pi s n}\right)^{1/2}\\
& =   \left( \cosh \frac{\pi s n}4 \right)^{1/2}   \le {e^{ \frac{\pi |s| n}8}}.
\end{align*}

Therefore, regardless of the parity of $n$ we have $ A(s) \le {e^{ \frac{\pi |s| n}8}}$.
Denoting  $A_0(s) = {e^{ \frac{\pi |s| n}8}}  $ and $A_1(s)=1$, we obtain
 $$\|I_{\frac{n}{2}is}f\|_{\infty} \le A_0(s) \|f\|_1$$
and 
$$\|I_{\frac{n}2(1+is)}f\|_{2} \le A_1(s) \|f\|_2$$
for all $s\in \mathbb R$.

The   Stein interpolation  theorem now implies
that for $0<p<n/2$ we have
\begin{equation}\label{bound}\|I_{p}f\|_{n/p} \le C(n,p) \|f\|_{n/(n-p)},  
\end{equation} 
for some constant $C(n,p)$.

Instead of using Stein's bound for $C(n,p)$, we will proceed as follows. 
Consider the function $F(z)=- \frac{n}4 \Re (z\log z)$, which is clearly harmonic in the strip $0< \Re(z)< 1$.
Now let us compute $F(z)$ on the boundary of the strip. When $z= is$, we get $\log(  is)=\log | s| + \frac{i\pi}{2}$ if $s>0$ and $\log( is)=\log |s| - \frac{i\pi}{2}$ if $s<0$. Therefore,
$$F(  is)=- \frac{n}4 \Re \left( is\left(\log | s| +\mbox{sgn}(s) \frac{i\pi}{2}\right)\right)= \frac{\pi n}{8}|s|.$$
Observe that $F(is)=\log A_0(s).$

When $z=1+is$, we get
\begin{align*}F(1+is)&=-\frac{n}4 \Re \left((1+is)\left(\log \sqrt{1+s^2} +i\arctan s \right)\right)\\ & = -\frac{n}8\log(1+s^2)+\frac{n}4 s \arctan s .
\end{align*}
The derivative of the latter function is negative when $s<0$ and positive when $s>0$. Therefore the function achieves its minimum at zero, and thus $F(1+is)\ge 0=\log A_1(s)$ for all real $s$.

Summarizing, on the boundary of the strip the following holds: $\log A_0(s)=F(is)$ and $\log A_1(s)\le F(1+is)$, for all $s\in \mathbb R$.  We can now conclude that for all $p\in (0,\frac{n}2)$ we have
$$\log C(n,p)\le F\left(\frac{2p}{n}\right)=-\frac{p}2 \log\left(\frac{2p}{n}\right).$$
That is, 
$$C(n,p)\le \left(\frac{n}{2p}\right)^{\frac{p}2},$$
which together with (\ref{bound}) yields the result.

\qed

We will now show that (\ref{MainIneq}) holds up to a multiplicative constant (depending on $n$ and $p$) for all origin-symmetric star bodies.

\begin{theorem}
For  every
$0< p < n/2$ and every origin-symmetric star body $K$ in $\mathbb R^n$ we have
	$$ \int_{S^{n-1}} \left| I_p (\|\cdot\|_K^{-n+p}) (\theta) \right|^{n/p} d\sigma(\theta) \le \left(\frac{n}{2p}\right)^{\frac{n}2}   \left(\kappa_n\right)^{1-\frac{n}{p}} (\mathrm{vol}_{n}(K))^{\frac{n}{p}-1}. 
	$$

\end{theorem}\label{isom}
\proof This is a direct application of Theorem \ref{interp} to the function $f=\|\cdot\|_K^{-n+p} $.	 
\begin{align*}& \left(  {\int_{S^{n-1}} \left| I_p (\|\cdot\|_K^{-n+p}) (\theta) \right|^{n/p}  d\sigma(\theta)}\right)^{p/n}\\
	&\qquad\le \left(\frac{n}{2p}\right)^{\frac{p}2}   \left( {\int_{S^{n-1}} \left( \|\theta\|_K^{-n+p}  \right)^{n/(n-p)}  d\sigma(\theta)}\right)^{(n-p)/n} \\
	&\qquad= \left(\frac{n}{2p}\right)^{\frac{p}2}   \left( \kappa_n^{-1} \mathrm{vol}_{n}(K) \right)^{(n-p)/n}  .
	\end{align*}
	Raising both sides to the power $n/p$, we get the result.

\qed 

\begin{remark} Observe that the result above improves inequality (\ref{KPZ})   when $k>
c\frac{n}{\log^2 n}$.    Indeed, let $0<k<n/2$ and let $K$ be an    origin-symmetric star body in $\mathbb R^n$ whose $k$-intersection body exists. If we assume that $\mathrm{vol}_{n}(K)=\mathrm{vol}_{n}( B_2^n)$, then Theorem \ref{isom} yields
$$ \left(\frac{\mathrm{vol}_{n}(I_k (K))}{ \mathrm{vol}_{n}(I_k(B_2^n))}\right)^{1/n}  \le \sqrt{\frac{n}{2k}} .$$
 Note that we do not require  $K$ or  $I_k (K)$ to be convex. 
 \end{remark}

Before solving a local version of Question \ref{Conj1} we will prove the following  lemma pertaining to the equality case in Question \ref{Conj1}. 

\begin{lemma} The conjectured inequality (\ref{MainIneq})   becomes equality if $K$ is a centered ellipsoid. 
\end{lemma} 
\proof
First, let us show that  inequality (\ref{MainIneq}) is invariant under invertible linear transformations.
Indeed, let $T\in GL_n(\mathbb R)$. Applying the transformation $T$ to the right-hand side of  inequality (\ref{MainIneq})    yields a factor of $|\det T|^{n/p-1} $. Now let us show that the same happens to the left-hand side of 
(\ref{MainIneq}).

 Consider the origin-symmetric star-shaped set $L$ defined by the formula
$$\|\theta\|_L^{-p} =  \left|  \left( \|\cdot\|_K^{-n+p} \right)^\wedge(\theta)\right|, \quad \theta \in S^{n-1}.$$ Using the connection between the Fourier transform and linear transformations, we get

\begin{align*}
& \int_{S^{n-1}} \left| (\|\cdot\|_{TK}^{-n+p})^\wedge(\theta) \right|^{n/p} d\theta = \int_{S^{n-1}} \left| (\|T^{-1}(\cdot)\|_{K}^{-n+p})^\wedge(\theta) \right|^{n/p} d\theta  \\
&\quad  = |\det T|^{n/p} \int_{S^{n-1}} \left| (\|\cdot \|_{K}^{-n+p})^\wedge(T^{t} \theta) \right|^{n/p} d\theta  \\
&\quad   = |\det T|^{n/p} \int_{S^{n-1}}  \|T^{t} \theta \|_{L}^{-n}  d\theta = |\det T|^{n/p} \int_{S^{n-1}}  \| \theta \|_{ T^{-t}  L}^{-n}  d\theta \\
&\quad   =|\det T|^{n/p}  n \mathrm{vol}_n(T^{-t} L)=|\det T|^{n/p-1}  n \mathrm{vol}_n( L)\\
&\quad   =|\det T|^{n/p-1} \int_{S^{n-1}} \left| (\|\cdot\|_{K}^{-n+p})^\wedge(\theta) \right|^{n/p} d\theta,
\end{align*}
where we used $T^t$ and $T^{-t}$ to denote the transpose and the inverse of the transpose of $T$ correspondingly. The calculations above remain valid if the Fourier transform is replaced by the operator $I_p$ since they are equal up to a constant multiple. Thus, the linear invariance of (\ref{MainIneq})  follows.

It remains to show that (\ref{MainIneq})   turns into equality when $K$ is a unit ball. Indeed,   we have
\begin{equation*}   \int_{S^{n-1}} \left| I_p(|\cdot|_{2}^{-n+p}) (\theta) \right|^{n/p} d\sigma(\theta)=1 =\left(\kappa_n\right)^{1-\frac{n}{p}} \left(\mathrm{vol}_n(B_2^n) \right)^{\frac{n}{p}-1}.
\end{equation*}

\qed

We will now prove that a local version of Question \ref{Conj1} has a positive answer.

\begin{theorem}
	Let $0<p<n/2$   and let  $K$ be an origin-symmetric star body in $\mathbb R^n$. If $K$ is sufficiently close to the Euclidean ball in the Banach-Mazur distance, then
\begin{equation}\label{C11} \int_{S^{n-1}} \left| I_p (\|\cdot\|_K^{-n+p}) (\theta) \right|^{n/p} d\sigma(\theta) \le     \left(\kappa_n\right)^{1-\frac{n}{p}} \left(\mathrm{vol}_n(K) \right)^{\frac{n}{p}-1},  
\end{equation}
with equality if and only if $K$ is an  ellipsoid centered at the origin.

\end{theorem}

\proof

Since $K$ is close to the Euclidean ball in  the Banach-Mazur distance, we can apply an appropriate linear transformation and assume that $K$ is close to $B_2^n$ in the Hausdorff distance, that is,
\begin{equation}\label{inclusion}
(1-\varepsilon) B_2^n \subset K\subset  (1+\varepsilon) B_2^n
\end{equation}
for some small $\varepsilon >0$.

Applying another linear transformation we can put $K$ in isotropic position, i.e., a position for which the following holds: 
$$ \int_K x_i x_j \, dx =  \lambda \delta_{ij},$$
for some positive constant $\lambda$ and all $1\le i,j\le n$; see \cite[Section  2.3.2]{BGVV} for details.

After putting $K$ in isotropic position, one can check that it is still close to $B_2^n$ in the Hausdorff distance; see \cite[Section 4]{ANRY}. So we can assume that (\ref{inclusion})   holds (with a different $\varepsilon$).

Let us write 
\begin{equation}\label{norm-n+p}
\|x\|_K^{-n+p} = H_0( 1+\varphi (x)), \quad x\in S^{n-1},
\end{equation}
 where $H_0$ is a constant (the harmonic of order zero in the spherical harmonic expansion of $\|x\|_K^{-n+p}$) and $\int_{S^{n-1}} \varphi(x) d\sigma(x)=0.$
 
Note that (\ref{inclusion}) implies 
  $$ (1+\varepsilon)^{-n+p} \le \|x\|_K^{-n+p} \le  (1-\varepsilon)^{-n+p}$$
  for all $x\in S^{n-1}$. Therefore, 
   $$ (1+\varepsilon)^{-n+p} \le H_0 \le  (1-\varepsilon)^{-n+p}.$$

 Since  $H_0$ is  close to one, if we dilate $K$ by a factor of $({H_0})^{1/(n-p)}$, $K$ will still be close to $B_2^n$ in the Hausdorff metric. So  from now on we will assume that  
    \begin{equation}\label{norm}
   \|x\|_K^{-n+p} = 1  + \varphi(x),
   \end{equation}
   where $\max_{x\in S^{n-1}}|\varphi(x)| < \varepsilon$ and $\int_{S^{n-1}} \varphi(x) d\sigma(x) =0.$

Let $$ \sum_{m\ge 2, \atop m\text{
		even}}   H_m$$ be the spherical harmonic expansion of $\varphi$.

Recall that $K$ is in isotropic position. We will show that this implies that 
\begin{equation}\label{H2}
\|H_2\|_2\le C\varepsilon \|\varphi\|_2. 
\end{equation}
Indeed, let $H(x)= \sum_{i,j=1}^n a_{ij} x_i x_j$ be a harmonic quadratic polynomial on $\mathbb R^n$. Observe that we necessarily have $\sum_{i =1}^n a_{ii}=0$.

Therefore,
\begin{multline*} \int_{S^{n-1}} \|x\|_K^{-n-2} H(x) \, dx = (n+2) \int_K H(x)\, dx    \\ =  (n+2)\sum_{i,j=1}^n a_{ij}  \int_K x_i x_j \, dx =  (n+2)\sum_{i,j=1}^n a_{ij} \lambda \delta_{ij} =0.
\end{multline*} 
Thus  $\|x\|_K^{-n-2}$ has no second order harmonic in its spherical harmonic expansion. 

Raising (\ref{norm}) to the power $(n+2)/(n-p)$ and using the Taylor expansion, we get $$\Big|\|x\|_K^{-n-2} - 1 - \frac{n+2}{n-p} \varphi(x)\Big|  \le C \varepsilon |\varphi(x)|.$$
Taking the $L^2$-norms of both sides and keeping only the second order harmonic in the left-hand side, we obtain (\ref{H2}).

We will now compute the left-hand side of (\ref{C11}).
To this end, applying $I_p$ to both sides of  (\ref{norm}) and raising to the power $n/p$, we get
\begin{equation}\label{auxi}
\left|I_p(\|\cdot\|_K^{-n+p}) (\theta)\right|^{n/p} = \left|1 +I_p\varphi(\theta)\right|^{n/p},
\end{equation}
for all $\theta\in S^{n-1}$.

To expand the right-hand side of the latter equality, we will need the following observation.
For a fixed $\alpha>2$  let $\zeta=\zeta_\alpha$ be  the function defined by $\zeta(t) = |1+t|^\alpha-1-\alpha t - \frac{\alpha(\alpha-1)}{2}t^2$, for $t\in \mathbb R$.
We claim that
\begin{equation}\label{zeta}\left|\zeta(t)\right|\le \left\{\begin{array}{ll} D  |t|^\alpha, & \mbox{if } 2<\alpha\le 3, \\  D  (|t|^3+|t|^\alpha), & \mbox{if } 3\le\alpha, \end{array}\right.
\end{equation} 
where $D$ is a constant (depending on $\alpha$). 
To prove the claim, observe that $\zeta(t)$  divided by the right-hand side of (\ref{zeta}) is a continuous function of $t\in \mathbb R\setminus\{0\}$ with finite limits when $t\to 0$ and $t\to \pm \infty$.

Thus,
\begin{multline*} 
\int_{S^{n-1}} \left|I_p(\|\cdot\|_K^{-n+p}) (\theta) \right|^{n/p} d\sigma(\theta)  = \int_{S^{n-1}} \Big| 1 +I_p\varphi(\theta) \Big|^{n/p} d\sigma(\theta)\\
 = \int_{S^{n-1}} \Big( 1 +\frac{n}{p} I_p\varphi(\theta)+\frac{n(n-p)}{2p^2} \left(I_p\varphi(\theta)\right)^2 +\zeta\left(I_p\varphi(\theta)\right)\Big)  d\sigma(\theta).
\end{multline*}

Since $\varphi$ has no spherical harmonic of degree zero, neither does $I_p\varphi$. That is,  
$$ \int_{S^{n-1}} I_p\varphi (\theta)d\sigma(\theta)  =0.$$
Denoting 
$$R_1=\int_{S^{n-1}}  \zeta\left(I_p\varphi(\theta)\right) d\sigma(\theta),$$
we get
\begin{equation}\label{LHS}\int_{S^{n-1}} \left|I_p(\|\cdot\|_K^{-n+p}) (\theta) \right|^{n/p} d\sigma(\theta)  = 1 +   \frac{n(n-p)}{2p^2} \|I_p\varphi\|_2^2 +R_1.
\end{equation}

We will now show that $R_1=o(\|{\varphi}  \|_2^2)$.  By Theorem \ref{interp} there is a constant $C$ (depending on $n$ and $p$) such that 
$$\|I_{p}\varphi\|_{n/p} \le  C  \|\varphi\|_{2}. $$
%  Let us denote $\alpha = p/n$
When $2< n/p \le 3$, by (\ref{zeta}) we have 
$$|R_1|\le D \int_{S^{n-1}} \left| I_p\varphi(\theta) \right|^{n/p} d\sigma(\theta) =  D \|I_p\varphi\|_{ n/p}^{ n/p}\le C^{n/p} D \|\varphi\|_{2}^{ n/p}=o(\|{\varphi}  \|_2^2).$$
When $3\le n/p$, (\ref{zeta}) yields
\begin{align*}|R_1|&\le D\left( \int_{S^{n-1}} \left| I_p\varphi(\theta) \right|^{n/p}d\sigma(\theta) +\int_{S^{n-1}} \left|I_p\varphi(\theta) \right|^{3} d\sigma(\theta)\right)\\
 &=D \left(\|I_p\varphi\|_{n/p}^{n/p}+\|I_p\varphi\|_{n/p}^{3}\right)\\
& \le D \left(C^{n/p} \|\varphi\|_{2}^{n/p}+ C^3\|\varphi\|_{2}^{3}\right)=o( \|\varphi\|_{2}^{2}) .
\end{align*}
Thus, in both cases, $R_1=o( \|\varphi\|_{2}^{2}).$

We will now compute the right-hand side of (\ref{C11}). Using (\ref{norm}) we get
$$ \|x\|_K^{-n} =  \left(1 +   \varphi(x )  \right)^{n/(n-p)} = 1 +  \frac{n}{n-p} \varphi(x) +   \frac12 \frac{np}{(n-p)^2}  \varphi^2(x) +\eta(x)    , $$
where 
$$|\eta|\le c   \varepsilon  \varphi^2  ,$$
for some constant $c$.

Using that the integral of $\varphi$ over the sphere vanishes, we get \begin{align*}\mathrm{vol}_n(K)& = \kappa_n \int_{S^{n-1}} \|x\|_K^{-n}\, d\sigma(x) \\ & =  \kappa_n\left(1+   \frac{n p}{2(n-p)^2}     \int_{S^{n-1}} \varphi^2(x)\,d\sigma(x) + R_2\right), 
\end{align*}
where 
$R_2=o( \|\varphi\|_{2}^{2}).$

Hence,
\begin{equation}\label{RHS}    \left(\kappa_n\right)^{1-\frac{n}{p}}\mathrm{vol}_n(K)^{\frac{n}{p}-1} =1 +  \frac{ n }{2(n-p) }  \|\varphi\|_{2}^{2} + R_2,
\end{equation}
where $R_2$ is different from that above, but is still of order
$o( \|\varphi\|_{2}^{2}).$

Let us now compare   (\ref{LHS}) and  (\ref{RHS}). Since   $R_1$ and $R_2$ are of order $o( \|\varphi\|_{2}^{2})$,
to finish the proof we need to show that
\begin{equation}\label{to-prove} \frac{n(n-p)}{2p^2} \|I_p\varphi\|_2^2  \le \frac{ n }{2(n-p) } \|\varphi\|_2^2  +o( \|\varphi\|_{2}^{2}),
\end{equation}
provided $\|\varphi\|_{2}$ is sufficiently small.

Indeed,
\begin{align*} \frac{n(n-p)}{2p^2} & \|I_p\varphi\|_2^2  = \frac{n(n-p)}{2p^2}\left(\lambda_2^2(n,p)  \|H_2\|_2^2+ \sum_{m\ge 4, \atop m\text{
		even}} \lambda^2_m(n,p) \|H_m\|^2_2\right)\\
	& \le \frac{n(n-p)}{2p^2}\left(\lambda_2^2(n,p)  \|H_2\|_2^2+ \lambda^2_4(n,p) \sum_{m\ge 4, \atop m\text{
			even}}  \|H_m\|^2_2\right)\\
			&  = \frac{n(n-p)}{2p^2}\left((\lambda_2^2(n,p) -\lambda^2_4(n,p) )  \|H_2\|_2^2+ \lambda^2_4(n,p) \sum_{m\ge 2, \atop m\text{
				even}}  \|H_m\|^2_2\right)\\
	& = o(\|\varphi\|_2^2)+ \frac{n(n-p)}{2p^2}   \lambda^2_4(n,p) \|\varphi\|_2^2\\
	&= o(\|\varphi\|_2^2)+ \frac{n}{2(n-p)}   \frac{(p+2)^2}{(n-p+2)^2}  \|\varphi\|_2^2 .
	\end{align*}
	Since $\frac{(p+2)^2}{(n-p+2)^2} <1$,    (\ref{to-prove}) follows.
	
	\qed

{\bf Acknowledgment.} The author is grateful to Fedor Nazarov and Dmitry Ryabogin for their invaluable help.

\end{document}